\documentclass{amsart}

  \usepackage[all]{xy}

  \usepackage{epsf,epsfig,amsfonts,graphicx,color}

\usepackage{amsmath,amssymb}
  \usepackage{url}
   
 \def\R{\mathbb R}

 \def \E {\mathbb E}

 \def \S {\mathbb S}

\newcommand{\beq}{\begin{equation}}
\newcommand{\eeq}{\end{equation}}
\newcommand{\Leq}[1]{\label{#1}\end{equation}}

\newtheorem{theorem}{Theorem}[section]
\newtheorem{proposition}{Proposition}[section]
\newtheorem{corollary}{Corollary}[section]
\newtheorem{lemma}{Lemma}[section]
\newtheorem{definition}{Definition}[section]

\newtheorem{conj}{Conjecture}

\begin{document}

\title{The negative energy  N-body problem has finite diameter}

\author{Richard Montgomery}
\address{Mathematics Department\\ University of California, Santa Cruz\\
Santa Cruz CA 95064}
\email{rmont@ucsc.edu} 

\maketitle

\begin{abstract} The Jacobi-Maupertuis metric provides a reformulation 
of the  classical N-body problem   as a  geodesic flow on an  energy-dependent
metric space denoted $M_E$ where  $E$ is the    energy of the problem. 
We show   $M_E$
has finite diameter for $E < 0$. Consequently  $M_E$
 has no ``metric rays''.  Motivation comes from   work of   Burgos- Maderna and Polimeni-Terracini  
 on the case $E \ge 0$ 
and  from a need to correct an error made in a   previous ``proof''. 
$M_E$ is constructed by completing the JM metric, a Riemannian metric on the Hill region, a domain in configuration space. 
We show that $M_E$  has finite diameter for $E < 0$ by showing that there is  
a constant $D$ such that all points of the Hill region  lie a distance $D$ from  the Hill boundary. (When $E \ge 0$ 
the Hill boundary is empty.)  
 The proof relies  on a game of escape  which allows us to   quantify the escape rate  from a closed subset
 of configuration space, and the reduction of this game to one of  escaping  the boundary of a polyhedral convex cone
  into its interior. 
  \end{abstract} 

\date{today}

The  classical N-body problem at  fixed energy  $E$ can be reformulated
as a geodesic problem.  The geodesics are those of  the   Jacobi-Maupertuis [JM] metric  at  that energy.
  Our main result is
\begin{theorem}  The Jacobi-Maupertuis [JM]metric for the N-body problem at negative energy has  finite diameter.
\label{thm 1}
\end{theorem}

The configuration space of the N-body problem is a Euclidean space endowed
with a potential energy 
function $q \mapsto V(q)$.   The Jacobi-Maupertuis
metric at energy $E$ is defined on the region   $\{q:   - \infty <   V(q) < E \}$ whose 
closure,   the domain   $\{q:  -\infty \le V(q) \le E \}$,  we  call the Hill region.
The boundary   $\{ q: V(q)  = E \}$ is called the Hill boundary and will be denoted as $\partial M_E$.  
 The Hill boundary is  non-empty if and only if $E < 0$.  The  Jacobi-Maupertuis
  metric degenerates to zero at the Hill boundary and so  we can travel along the Hill boundary at zero (metric) cost.   
Write $M_E$ for the metric completion of the Jacobi-Maupertuis metric
and $d_E$ for the corresponding metric.  When forming $M_E$ for $E <0$ 
we must  collapse the Hill
boundary to a single point because of this zero-cost travel.  We also  denote this point  by  $\partial M_E$.     
 The theorem above follows   from
  \begin{theorem}  If $E < 0$ then $sup_{\{q \in M_E\}}  d_{E} (\partial M_E,  q) < \infty$
 \label{thm: radius} 
 \end{theorem}
  
{\sc Proof of theorem \ref{thm 1} from theorem \ref{thm: radius}.}
The theorem asserts that  the boundedness of the function
$d_{E} (\partial M_E,  q)$. In other words,  there is a positive constant $K$ such
for all $q \in M_E$ we have that  $d_{E} (\partial M_E,  q) < K$.
Take  $p, q \in M_E$.  Then  
$
d_E (p,q) \le d_E (p, \partial M_E) + d_E (\partial M_E, q) \le 2K
$
so that the diameter of $M_E$ is less than or equal to $2K$. 

\vskip .3cm

{\sc remark 1.}  $M_E$ is well-known to have  infinite diameter when $E \ge 0$.    
 
 \subsection{Motivation}

The JM metric perspective on N-body dynamics has recently proven  particularly  useful
for  positive energies.  Maderna and Venturelli  \cite{Maderna1}  combined  this  perspective
  with   weak KAM methods    to prove  that given any  
initial configuration and any final asymptotic `hyperbolic' state  that there is a 
positive energy solution  connecting the two.   
%The completion of the Hill region with its JM metric will be denoted $M_E$. 

Maderna and Venturelli's positive energy solutions are  metric rays in $M_E$, $E > 0$.  
A {\it metric ray}  in a metric space $M$ is an isometric embedding
of the  half-line $[0, \infty) \subset \R$  into that metric space, the half-line being endowed with the
usual metric   inherited from $\R$.   A minimizing
geodesic  $c:[a, b] \to M$ is an isometric embedding of the interval $[a,b]$ into $M$.
If  $\gamma:[0, \infty) \to M$ is a metric ray then the restriction of $\gamma$
to any compact subinterval $[a,b] \subset [0, \infty)$   is a minimizing geodesic.   
(See Burago et al \cite{Burago} for general definitions and results concerning geodesics in metric spaces.) 

Burgos and Maderna (\cite{Burgos2}) obtained a  `parabolic generalization'' of Maderna-Venturelli,
and in so doing asked the compelling  question: {\bf  If the energy $E$ is negative  does
$M_E$ admit any metric rays?}  
 %(See also Polimeni-Terracini \cite{Terracini} regarding
%JM approaches to parabolic-hyperbolic escapers.)  
Since finite diameter spaces cannot support metric rays
 our theorem provides   an immediate answer to their question. 
   \begin{corollary}  There are no metric rays in $M_E$ for $E < 0$.  
 \end{corollary}
 
 \begin{figure}
  \includegraphics[width=8cm]{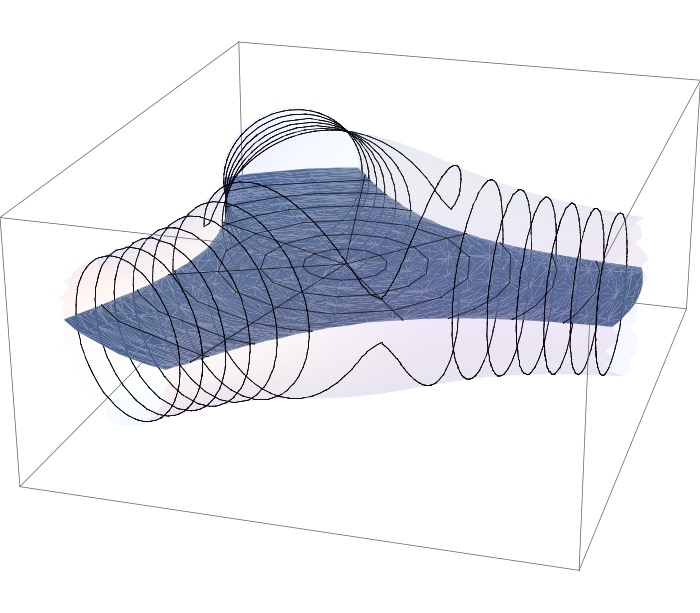}
 \caption{The Contour level surface  $V = -1$   drawn in the planar 3-body shape space
 is like the  surface of a  plumbing fixture consisting of three pipes  centered
 about the three  binary collision rays. 
 The Hill region   $\{ V \le -1 \} $
 projects onto the interior of the surface.    The  shaded planar domain  inside the Hill region is   the 
  Hill region for the collinear three-body problem at energy $-1$.  (Courtesy of Rick Moeckel.)} 
 \label{fig: HillRegion2}
 \end{figure}

      On page 387 of \cite{MontRCD} I had claimed to have proven     theorems \ref{thm 1} and
\ref{thm: radius}.  
  That proof is wrong.   I also  wrote this article  to correct that error.  
  In detail, the error is as follows.  I gave an estimate
  showing that   any point could be
connected to the Hill boundary by a Euclidean ray whose Jacobi-Maupertuis length
was   linear in $\|A \|$
where   $A$ is  the point where that ray pierces the   Hill boundary.  But  
 $\|A \|$ is unbounded on the Hill boundary,
so the claimed proof   does not bound $d_E( q, \partial M_E)$.

 \subsection{A conjecture}

 \begin{conj}   $sup_{\{q \in M_E\}}  d_{E} (\partial M_E,  q) =   d_E  (\partial M_E, 0)$ when $E < 0$.
 \end{conj} 
 This conjecture holds for $N=2$ bodies.  We can also    show that it holds
 locally for the 3-body problem.  
 `Locally' here means that  $q =0$ is a local maximum for $d_E (q, \partial M_E)$,
 and that $d_E (q, \partial M_E) < d_E (0, \partial M_E)$
for all  $\|q \|$ is sufficiently large.   
 $d_E  (\partial M_E, 0)$,  can be computed explicitly in
 terms of `minimal central configurations'' in a manner similar to the way that
 the minimal action to total collision in a fixed time can be computed by `dropping'
 minimal central configurations.   See for example \cite{Brake-to-Syz}.

 \section{Set-up.  Jacobi-Maupertuis metrics and the N-body equations.}

Let $\E$ be a real  inner product space and $V$ a smooth real-valued function on $\E$.   Together this data defines  a Newton's equations:  
\beq
\ddot q = - \nabla V (q).   
\Leq{Newtons}
with conserved energy   
 $$E(q, v)  = K (v) + V(q),  \text{ where } \qquad K(v) = \frac{1}{2} \langle v, v \rangle  ,  v = \dot q .$$
 Here the dots denote time derivatives.
 The inner product $\langle v, v \rangle$ used to define the kinetic energy $K$  is the given inner product on $\E$.
The gradient ``$\nabla$'' in Newton's equations  is   relative to this inner product
 so that $dV(q)(h) = \langle \nabla V (q) , h \rangle$.   
 
 The function $V$ is called the potential energy.  In order to accomodate the N-body problem we allow for points $q$ with 
  $V(q) = - \infty$.  We call these collision points.  The  gradient of $V$    goes to   infinity as we approach   
a collision point
 so   Newton's equations break down. 
 
Since $K \ge 0$ we have
 that $V(q) \le E$ if a solution $q(t)$ has energy $E$.  Define
 the Hill region at energy $E$ to be the locus
 $$\text{ Hill region}  = \{q \in \E:  V(q) \le E \}$$
 and its boundary, called the Hill boundary,  to be
 $$\partial M_ E : = \{q \in \E:  V(q)  = E \}$$
 All energy $E$ solutions $q(t)$ to Newton's equations must lie in the Hill region.  If the  solution $q(t)$  encounters the
 Hill boundary at some instant $t = t_0$  then $K(\dot q (t_0)) = 0$  which means the solution has instantaneously
 stopped.    We call such instants or locations along the path ``brake points''. 
 
   The Jacobi-Maupertuis principle
 asserts that the solutions to Newton's equations at energy $E$ can be characterized
 as geodesics on the Hill region.

 \begin{definition} The Jacobi-Maupertuis  metric at energy  $E$
 for the Newton's equation (\ref{Newtons})   is the  Riemannian metric 
 \beq
 ds^2 _{JM} = 2(E-V(q)) \langle dq, dq \rangle 
 \Leq{JM metric} defined on the  interior  $\{-\infty <  V  <  E \}$ of the Hill region for that energy,
 with collisions ($V = - \infty$) excluded.   
  \end{definition} 
  
  Note that the conformal factor $E-V(q)$ vanishes on the Hill boundary: the metric fails to be Riemannian
  and paths can travel for no cost along the Hill boundary.  
  
 \begin{theorem} [Jacobi-Maupertuis principle]  Away from   collisions  and brake  points, the energy $E$  solutions to Newton's equations
are reparameterizations of  geodesics on the    Hill region.     Conversely, away from the collisions  and   the Hill boundary,   geodesic
for this Jacobi-Maupertuis metric are reparameterizations of 
 energy $E$ solutions of Newton's equations.
 \end{theorem} 
 
 \noindent See 
 Landau-Lifshitz \cite{Landau} p. 141,   Knauf \cite{Knauf}  p. 178,  or Abraham-Marsden \cite{AbMa} p.228  for proofs of the Jacob-Maupertuis principle.

 The metric completion  of the interior of the Hill region will be denoted as $M_E$.
 If the Hill boundary is non-empty and connected it contains that boundary, collapsed to a point,
 as a single point denoted as $\partial M_E$.  That point is typically not a manifold point.
 Depending on the potential $V$, the collision locus may, or may not be in $M_E$.  
 For a  detailed description  of $M_E$ see
 Proposition 1 on p. 383 of \cite{MontRCD}. 
 % XXX \rmont{ fill in REF} 
 
 \subsection{The N-body equations}

 We  put the N-body problem into the above framework by using the  set-up   which   Albouy and Chenciner 
 taught me. See \cite{Albouy-Chenciner} or \cite{Chenciner}.  View the bodies as point masses.
 They move in   $d$-dimensional Euclidean space $\R^d$.
 The standard choice of $d$ is   $d =3$.   The bodies, or masses, are labelled by an index   $a \in [N] = \{1, 2, \ldots, N\}$.
The  instantaneous location   of 
 the $a$th  body is denoted by 
 $q_a \in \R^d$ so that the  simulataneous positions of all N-bodies
 is    encoded by the  vector
 $$q = (q_1, q_2, \ldots, q_N) \in \E: = (\R^d )^N.$$
 The mass of the $a$th body is $m_a > 0$ and the masses endow   $\E$ with an inner product $\langle  \cdot  ,  \cdot  \rangle$
 which we call the mass inner product and 
whose associated quadratic form is :   
$$\langle  q ,  q \rangle =  \Sigma m_a |q_a |^2, $$
where $|q_a|^2 = q_a \cdot q_a$ is  the standard dot product on
$\R^d$.  
  When applied to velocities $\dot q = (\dot q_1, \ldots \dot q_N);  \dot q_a := dq_a/dt$
the mass inner product  yields twice the kinetic energy:
\begin{equation}
\begin{aligned} 
K(\dot q) & = \frac{1}{2} \langle \dot q , \dot q \rangle  \\
  & =  \frac{1}{2} \Sigma m_a | \dot q_a | 
  \end{aligned}
\end{equation}  
The potential energy $V = - U$ is the negative of
$$U(q) = G \Sigma \frac{ m_a m_b}{r_{ab}};  r_{ab} = |q_a - q_b|,  \qquad  V = - U$$
the sum being over all distinct pairs $a, b$ taken from $[N] = \{1, 2, \ldots, N\}$,
and $G$ being the gravitational constant.  The total
energy is  then given by: 
 \begin{equation}
\begin{aligned} 
E(q, \dot q) &  = K(\dot q)  - U (q)  \\
      \end{aligned}  
\label{energy}
\end{equation}

Since $E \ge V \iff U \ge -E$ we have that  the   Hill region for energy $E$
 is the domain $\{q :  U(q) \ge - E \}$ within $\E$.  
 Since  $U(q) > 0$ everywhere we have that
 the Hill region is all of $\E$ whenever $E  \ge 0$.
 On the other hand, if $E < 0$ the Hill region is not all of
 $\E$ and has a non-empty boundary
 $$\partial M_E = \{q:  U(q) = -E \}$$
 The JM metric on the interior of $M_E$ is  
 $$ds^2 _E =  2 (U(q) - E)  \langle dq, dq \rangle. $$ 
 
 We note that the metric is conformal to the flat Euclidean metric  $\langle dq, dq \rangle$
 but that the conformal factor is zero along the Hill boundary.
 It follows that we can {\it travel for free} along the Hill boundary.
 It costs nothing to move along the Hill boundary.

\subsection{Scaling Normalization}

We    reduce  to the case  of  energy  $E =-1$   using
the standard scaling symmetry  for the N-body problem.   This  scaling asserts that
if  the curve $q(t) \in \E$ is a solution to (\ref{Newtons}) then, for any positive real number $\lambda$,  so is  $\lambda q (\lambda ^{-3/2} t)$
and that if the 1st solution has energy $E$ then the second solution has energy  $E/\lambda$.  At the metric level, this 
scaling symmetry corresponds to
the fact that the scaling substitution
that $q = \lambda Q$ takes the JM metric on $M_{E/\lambda}$ to $\lambda^{1/2}$ times
the JM metric on  $M_E$.  By appropriate $\lambda$ we can scale any negative $E$ to $E =-1$. 
So, from now on, set $E=-1$, write $M_{-1} = M$
 and  the JM metric  $d_{-1}$ on $M$ as $d$.

\subsection{Collisions}  The collision locus $\Delta  \subset \E$ 
  is the union 
  $$\Delta = \bigcup_{ \text{distinct pairs} } \Delta_{ab}$$
  of the linear subspaces $\Delta_{ab} = \{r_{ab} = 0 \} = \{q \in \E:  q_a = q_b \}$.
 The potential blows up exactly at the points of  $\Delta$ and the forces, or gradients , the right hand side of Newton's equations,
 also blow  up along the collision locus.

\subsection{Distance to collision} 

The following lemma is crucial to our proof. 
\begin{lemma}  The Hill region $\{q:  U(q) \ge 1\}$ 
  lies a bounded Euclidean distance  $dist(q, \Delta)$ from the collision locus:   
     there is a positive constant $k$ such
that $U(q) \ge 1 \implies dist(q, \Delta) \le k$. 
\end{lemma}

 Here we write $dist(q, K) = inf_s (\|q - s \| : s \in K)$ 
 for $K$ a closed subset of $\E$.

The lemma   uses the fact that the function   $1/U(q)$ is Lipshitz equivalent to   
the function $dist(q, \Delta)$.  In other words, 
there exist positive constants $c_1   < C$ such that
that  for all $q \in \E$ we have 
\beq
c_1 dist(q, \Delta) \le \frac{1}{U(q)} \le C dist(q, \Delta).
 \Leq{ineq: dist potential}
 
 {\sc Proof of lemma.} The lemma follows immediately from the first inequality  of inequality
 (\ref{ineq: dist potential}) with $k = 1/c_1$.   QED

Inequality (\ref{ineq: dist potential})   is based on the distance formula \footnote{This equality 
is stated without a derivation in   the last paragraphs of my book on subRiemannian geometry.  For completeness,
I derive the  distance formula in an Appendix.}
\beq
dist(q, \Delta_{ab}) = k_{ab} r_{ab}, \text{ where }  k_{ab} = \sqrt{ m_a m_b/ (m_a + m_b)}
\Leq{distance fmla} 
from which it follows that
$$U = \Sigma \frac{ \lambda_{ab}} {dist(q, \Delta_{ab})} \qquad \text{ where }  \lambda_{ab} = G m_a m_b k_{ab} $$
Use
$$dist(q, \Delta) = min_{ab} dist(q, \Delta_{ab})$$
   to   get
$$\frac{\lambda_{ab}}{dist(q, \Delta_{ab})} < U(q) \le \Sigma \frac{ \lambda_{ab}} {dist(q, \Delta)}$$
valid for any pair $a, b$.  Choose a pair $a, b$ such that $dist(q, \Delta_{ab}) = dist(q, \Delta)$ 
and set 
$$\lambda_* = min_{a \ne b} \lambda_{ab} , \qquad \text{ and }  \Lambda = \sum \lambda_{ab}.$$
We get 
\beq
\frac{\lambda_*} {dist(q, \Delta)} \le U(q) \le \frac{\Lambda}{dist(q, \Delta)}
\Leq{dist ineq 2} 
which yields inequality (\ref{ineq: dist potential}) with constants  $c_1 = \frac{1}{\Lambda}$ and $C = \frac{1}{\lambda_*}$.

\section{A game of escape} 

Our proof of theorem \ref{thm: radius}  relies  on   a game of escape. 

{\bf Setting up the game.}  Select a  finite-dimensional Euclidean  space $\E$ together
with a finite collection $L_1, L_2, \ldots, L_k$ of distinct linear subspaces of $\E$,
no one of which is contained in any other.  

{\bf The game.}  The game is to  escape   the union of the $L_i$ as quickly as possible.

Write $\Delta = \cup_{i = 1} ^k  L_i$.
We measure escape in terms of  $dist(q, \Delta) = min_i dist(q, L_i)$.  
We use it to   quantify  the pay-off in the game which 
we call  the  ``escape rate''. 
We have found it helpful   to generalize the setting.

Let $\Delta$ be a closed   subset   of  Euclidean space $\E$.
For   $t > 0$ let 
$$N_t : = N_t (\Delta) := \{q \in \E:  dist(q, \Delta) \le t \}$$ denote   the set of points  of $\E$ lying within a distance $t$ of $C$.
Note that 
$$\partial N_t  = \{q :  dist(q, \Delta) = t \}.$$

{\sc Escape routes and their rates} 
\begin{definition}
A $t$ -escaper  is a rectifiable  path starting inside   $N_t(\Delta)$ and exiting $N_t (\Delta)$
 and along which the distance
from $\Delta$  is strictly monotonic increasing as a function of arc length. 
%The escape length $\ell$ of 
%a t-escaper  is  its length  up until first hitting $\partial N_t$.
\end{definition}

{\sc remark}  In order for the  distance from $\Delta$  to be strictly monotonic increasing
it must be true that   $\Delta$ has empty interior.

\vskip .3cm
  
We want to talk about the escape rate of an escaper.
\begin{definition} 
Parameterize a t-escaper $\gamma$ by   arclength $s$. Suppose that 
\beq
dist(\gamma(s), \Delta) \ge dist(\gamma(0), \Delta) + cs
\Leq{escape rate} 
holds for some   constant $c > 0$ and all   
 $s$ up to the escape time.  Then we will say that $\gamma$ has
t-escape rate at least  $c$.   The largest such $c$ will be called
the escape rate of the escaper $\gamma$. 
\end{definition}
\noindent And we want to talk about the escape rate from $\Delta$
\begin{definition}  If there is a   positive number $c$
such that  for all $t > 0$ and   all  $p \in N_t$ 
there exists  a t-escaper starting  at $p$ with  escape rate at least $c$ 
then we will say that the escape rate from $\Delta$ is positive and at least $c$.  The supremum of all
such  $c$'s is the escape rate from $\Delta$.   When we want to  distinguish this escape rate from the earlier 
escape rates we have defined we will refer to it as the {\sc global escape rate.} 
\end{definition}

Our  goal in  playing the  game and  making all  these  definitions is to  prove:  
\begin{theorem}  If $\Delta$ is the union of a finite collection of proper linear subspaces
of a Euclidean vector space  then
the global escape rate from $\Delta$ is positive  and at least  
$1/ dist(0, \partial N_1 (\Delta))$.   The t-escapers can be taken to be line segments. 
\label{thm:  subspace arrangements}
\end{theorem}
\noindent We postpone the  proof.

\section{Proof of theorem \ref{thm: radius} from theorem \ref{thm: subspace arrangements}}

Take $\Delta$ to be the collision locus.  
Inequality (\ref{ineq: dist potential}) shows that  if $dist(q, \Delta) \ge \frac{1}{C}$ then  $U(q) \le 1$.
In other words, a $t$-escaper  with  $t = 1/C$ has left the Hill region $\{ U \ge 1 \}$ and
 so has crossed the Hill boundary $\{ U = 1 \}$ at or before escape from $N_{1/C} (\Delta)$.
  Theorem \ref{thm: subspace arrangements}  applied to $\Delta$ 
 with $t = 1/C$  guarantees  a positive  global t-escape rate of  $c >0$.
   Now $U > 1$ implies  $dist(q, \Delta) < \frac{1}{C}$ :  the Hill region is contained in $N_{1/C} (\Delta)$.
Theorem \ref{thm: subspace arrangements} guarantees that  through any point $q_*$ of the Hill region  
there is a unit speed line segment of  $q(t)$ starting at $q_*$, satisfying  
$dist(q(t),  \Delta) \ge dist(q_*, \Delta) + ct$ and crossing the Hill boundary.
 Inequality (\ref{ineq: dist potential})
 now implies that $\frac{1}{U(q(t))}  \ge  c_1 (\frac{1}{U(q_*)} + ct)  \ge k t$
with  $k =  c_1 c$.   Thus
$$U(q(t)) \le \frac{1}{kt}.$$
It follows that along our t-escapers $q(t)$ we have  $\sqrt{U(q(t)) - 1}   \le \sqrt{\frac{1}{kt} -1}$
whenever $U(q(t) \ge 1$.  For simplicity,
set $$\lambda(x) = \begin{cases} \sqrt{x-1},  x \ge 1\\
0,  x \le 1
\end{cases}.
$$
The   JM arclength  integrand can then be written $\lambda(U(q(t)) \|\dot q (t) \| dt$
and the inequality just establishes asserts that 
%\rmont{ do I use `JM' before? change all instances of `JM' to long form?}
$\lambda(U(q(t)) \le \lambda({\frac{1}{kt}})$.
(Our escapers  $q(t)$ have escaped beyond the Hill boundary
at or before the time $t_* = 1/k$.)

Along an escaper   $\| \dot q (t) \| =1$
so that its   JM arclength satisfies 
$$\ell(q) \le \int_0 ^{1/k} \sqrt{ {\frac{1}{kt} -1} } dt  = k \int_0 ^1 \sqrt{ \frac{1}{u} -1} du  $$
This last  integral   is finite  and less that $1/2 = \int_0 ^1 \sqrt{ \frac{1}{u}} du$
so that the JM arclength to escape is less than $\frac{k}{2}$.
We have shown that any point inside the Hill region can be connected to the Hill boundary
by a path whose length is less than $k/2$, showing that
$dist (q, \partial M) \le c/2$ for all $q_*$.

QED

\section{Proof of the escape theorem, theorem  \ref{thm:  subspace arrangements}} 

%The proof of this simple-seeming theorem turned out to be a rather long story.  
We begin with  examples of  the escape game
and  consequent   escape rates.
 
\subsection{ Examples} \label{subsec: examples}
$$ $$

1. If $\Delta$ is a linear subspace then the escape rate from $\Delta$ is $1$. 
Lines orthogonal to $\Delta$ supply optimal escape routes.     
 
2. If $\Delta$ is the union of the x and y axis in the plane then 
$N_t (\Delta)$ is the union of the two strips $|x| \le t$ and $|y| \le t$.  The points $(t,t) \in \partial N_t (\Delta)$ and its mirrors
$(-t, t), (t, -t),  (-t, -t)$ in the   other three quadrants
are the furthest exit points from $\Delta$,  being the  furthest points on $\partial N_t (\Delta)$  from the origin.  The escape strategy 
  used in  the proof  of  of the polyhedral escape theorem \ref{polyhedral pf}  consists of a separate strategy
  for each quadrant.   In the first quadrant escape by translating  that quadrant toward the point $(t,t)$ at unit speed.
Thus $p$  in the first quadrant escapes along the path $p + s (1/\sqrt{2}, 1/\sqrt{2})$
The escape rate of these paths is $1/\sqrt{2}$.

3. If $\Delta$ is    the union of two lines  in the plane which make an angle of $\theta_* < \pi/2$ relative to each other
in the plane, then the escape rate  into the acute sector $0 \le \theta \le \theta_*$    is $\sin(\theta_* /2)$.   The strategy 
within that sector is the one of   the proof  of  of the polyhedral escape theorem \ref{polyhedral pf}: move according to
the translation $p+ s v$ where $v$ is the unit vector along the angle bisector  of this acute sector. 

4. Suppose that  $\Delta$ is the coordinate orthant in $\R^n$, by which we mean the union of the coordinate hyperplanes
$x_1 = 0$, $x_2 = 0$, etc. Then the escape rate from $\Delta$ is $1/\sqrt{n}$.  For escaping from the positive
orthant $x_i \ge 0$ mone at unit speed in the direction parallel to the ray  $x_1 = x_2 = \ldots = x_n$.  

\subsection{Wrong Turns and Gradient Flows}  
Originally, I had approached the problem  of trying to ``win'' the escape game by using the gradient flow of the function $F(q) = +dist(q, \Delta)$.
The trajectories of this flow, appropriately interpreted,
are piecewise linear t-escapers.  However, an example involving configurations near double binary collisions in the planar four-body
problem shows that the escape rates of the resulting solutions to $\dot q = \nabla F(q)$ can
be arbitrarily small, and so  in the end this approach to escaping was   not of help. 

A successful approach which yields the correct global escape rate  
is   based on  the gradient flow for $F_t (q) = -dist(q, \partial N_t (\Delta))$.
It uses the reduction to convex sets  described below
and the fact that the minus gradient flow to a convex set with non-empty interior
has trajectories which are shortest line segments to that convex set.  
The strategy we describe here (see the proof of the polyhedral escape theorem
\ref{polyhedral pf})  ended up being  significantly simpler than this one,
although its escapers   are the limits $t \to \infty$ of the t-escapers  of this alternative successful gradient flow approach.

\subsection{Reduction to Hyperplane arrangements}
We     reduce the proof of theorem   \ref{thm:  subspace arrangements} to the   case where the linear subspaces  $L_i$ 
of the theorem are  hyperplanes.
In order to do this, observe that
if $\Delta \subset \Delta'$ are closed subsets then $N_t (\Delta) \subset N_t (\Delta')$
and $dist(q, \Delta) \ge dist(q, \Delta')$.  
It   follows that if $\gamma$ is a t-escape path for $\Delta'$ 
which starts  in  $N_t (\Delta) \subset N_t (\Delta')$   then it is also
a   t-escape path for $\Delta$.  Consequently, if the escape rate from   $\Delta'$ is 
  $c' > 0$
then the  escape rate from $\Delta$ is also positive and   is some $c \ge c'$.  

Let $L_1, \ldots, L_i , \ldots, L_k$ be the linear subspaces.
Choose hyperplanes $H_i \supset L_i$
and write $\Delta' = \bigcup_i H_i$.
Then $\Delta \subset \Delta'$
so that $N_t (\Delta) \subset N_t (\Delta')$. 
It follows from the definitions that if $\Delta'$
has  global escape rate $c$ then the global escape rate from $\Delta$
is at least $c$. 

Relabel  the linear hyperplanes
to call them $L_i$ instead of $H_i$
and their union to be   $\Delta$. 
We have reduced our escape problem to: 

\begin{proposition} [Escaping linear hyperplane arrangements] If $\Delta = \bigcup_{i =1 }  ^k H_i$
is a hyperplane arrangement in a Euclidean vector space then
the escape rate from $\Delta$ is finite and greater than or each to $c_*$
where $c_*$ is given immediately following 
theorem  \ref{prop: polyhedral escape} below. 
\label{thm: hyperplanes} 
\end{proposition}

\subsection{Reduction to a Convex polyhedral cone.}

We  further   reduce     to the problem of escaping from a convex polyhedral cone.
Write each  hyperplane $L_i$ of proposition \ref{thm: hyperplanes}    in the form $L_i = ker(\ell_i)$ where 
$\ell_i \in \E^*$ is a nonzero linear functional.
Normalize $\ell_i$ to have unit length, so that 
$$\ell_i (q) = \langle q, n_i \rangle$$
  is the operation of inner product with respect to the unit normal vector $n_i$
  to $L_i$.   
%\rmont{ change H's back to L's since in future with K's use L's} 
%to $L_i$.   (There are two.)  
The distance of a point $q$ from $L_i$ is  $| \ell_i (q ) |$ so that 
$$dist(q, \Delta) = min_{ i  \in [k] } |\ell_i (q) |$$
where we use the symbol $[k]$ to mean the set $\{1, 2, 3, \ldots, k \}$.

 The complement of $\Delta$ consists
 of a finite number of  connected components:
 $$ \E \setminus \Delta = C_1 \cup C_2 \cup \ldots \cup C_M.$$
   Each component $C_{\alpha}$, $\alpha =1, 2, \ldots, M$  is an  open polyhedral cone
   bounded by some subcollection $L_i, i \in I = I(\alpha)$ of the hyperplanes making up our arrangement. 
 Any $t$-escaper must enter into one or another of these components.  Once entering it can never leave that component.     
 This `no exit' property  follows from   the fact that if  a continous curve is travelling   in a  component $C_{\alpha}$
  and then  exits to travel into   another component it must have crossed  one of the  bounding hyperplanes $L_i$.
  Crossing requires   $dist(q, \Delta) = 0$ which violates the assumption that  this distance function is strictly   monotone increasing along escapers.
  So, we have reduced to escapers which escape into the interior of a  single convex polyhedral cone. 
 
 In order to  focus on a single convex cone we  modify the definition of
 ``escape rate'' 
  to  focus on those paths  escaping into the cone's interior.
 \begin{definition}  Let $K$ be a convex cone with nonempty interior.
 Then by the escape rate into  $K$ we mean the global escape rate from $\Delta = \partial K$,
for  escapers which escape entering   the interior of $K$.
 \end{definition}
 
  \begin{theorem}
[Escaping polyhedral cones] Let $K$ be  a closed polyhedral cone with 
non-empty interior. 
Then the escape rate  into $K$  is positive.  The exact escape rate  is  $c= 1/ dist(0, K_1)$ where
 $K_1 \subset K$ is the convex polyhedron
defined by  $dist(q, \partial K) \ge 1$ for  $q \in K$. 
\label{prop: polyhedral escape}
\end{theorem}
\noindent

{\sc  Proof of 
Proposition \ref{thm: hyperplanes} }The escape rate from the union $\Delta$ of
the hyperplanes equals  the minimum of the escape rates into the interiors of 
all the   component  closed convex polyhedra $K^{\alpha}$ whose  interiors comprise
the components of  
the complement of $\Delta$.   In other words, for  each $K^{\alpha}$ take its  associated escape rate $c_{\alpha} = 1/ dist(0, (K^{\alpha} _1)$.
Let $c_* = min_{\alpha} c_{\alpha}$.  Since  the union of the $\partial K^{\alpha}$
form $\Delta = \bigcup_i L_i$, and since any t-escaper must escape into the interior of
one or another of these components, we have that the escape rate from
$\Delta$ is at least $c_*$.   QED

\subsection{Proving the polyhedral escape theorem}

We begin with some  background and notation.
A convex polyhedral cone $K \subset \E$.
  is  defined by a finite collection of linear inequalities:
 \beq
 K = \{ q:  \ell_i (q) \ge 0 , i  = 1, \ldots, m   \} 
 \Leq{Polyh cone}
 where the   $\ell_i \in \E^*$   are unit length linear covectors.
   Thus $\ell_i (q) = \langle n_i , q \rangle$
 where the  $n_i$ are  the inward-pointing unit normal vectors to   the cone's faces.
 Write 
 $$L_i =  \{ \ell_i = 0 \}$$
 for the corresponding hyperplanes
 so that  $F_i = L_i   \cap K$.
 The   distance of a point $q \in K$
 from a face $F_i$ is $\ell_i (q)$.
 Then , if $q \in K$ we have that 
 \beq
 dist(q, \partial K) =  min_{ i \in [m]} \ell_i (q), 
 % [m] = \{1, 2 , \ldots, m \}
 \Leq{eq: dist Delta}
where we are using the symbol $[m]$ for the set
of integers $1, 2, \ldots , m$.
  
  Define 
 \beq
 K_t =    \bigcap_{i \in [m]} \{ q:  \ell_i (q) \ge  t \}. 
 \Leq{star}
 $K_t$ itself is a convex polyhedron, being the 
  intersection of the  finite collection of half-spaces $\ell_i \ge t$.
 Observe that $K_t = \{ q \in K:  dist(q, \Delta) \ge t \}$ and that
 $\partial N_t = \partial K_t$ while $K \setminus int(N_t) = K_t$.

{\sc Remark.}  $K_1$, and hence $K_t$,  $t > 0$ is typically not a cone.
See Appendix B and especially accompanying  figure \ref{fig: coneEg}.    
\begin{figure}
  \includegraphics[width=8cm]{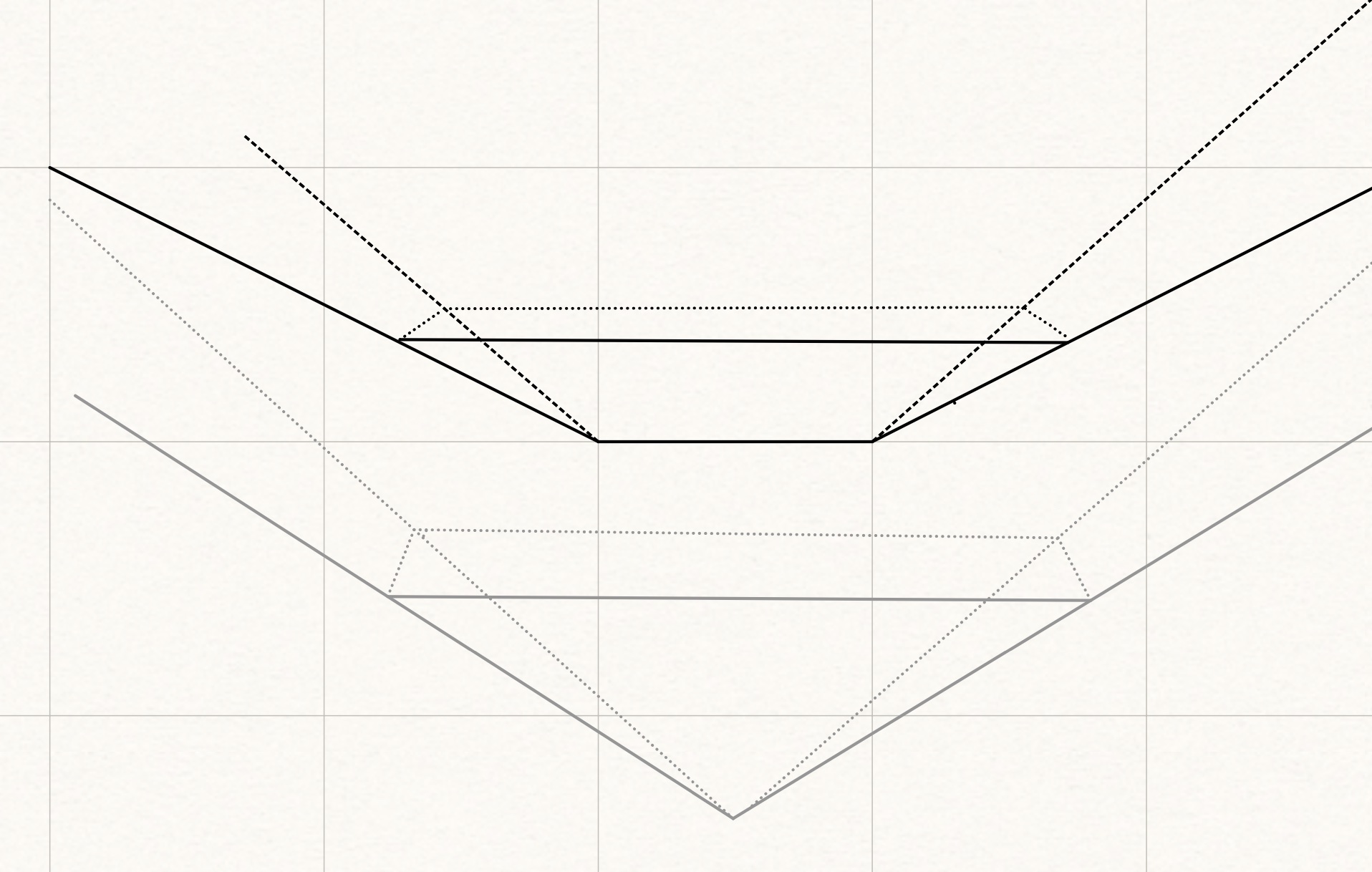}
 \caption{Cones $K_0$ (grey) and  its equidistant $K_1$ (black)  which is not a cone.
 See Appendix B for details.} 
 \label{fig: coneEg}
 \end{figure}
 
\subsubsection{ The proof of  polyhedral escape (theorem \ref{prop: polyhedral escape})} 
 \label{polyhedral pf}
 
In what follows $\Delta = \partial K$.       We will show that by translating $K$ inward we
obtain  a family of t-escapers covering $N_t (\Delta)$.

Let $\S$ denote the unit sphere in $\E$.   Choose any  unit vector $v \in \S \cap int(K)$.
Consider the one-parameter family of translations
$$p \mapsto \tau_s (p) = p + s v.$$   By convexity,  $\tau_s$  maps $K$ into $K$ for $s > 0$.
Write  $$c_* (v) = min_{i \in [m]} \ell_i (v) $$ 
I claim that for all $t >0$, the   family of rays $s \mapsto p + s v, s \ge 0$
 forms $t$-escapers with escape rate $c_* (v)$.    
 Indeed,
 $$\ell_j (p + s v) = \ell_j (p) + s \ell_j (v) \ge dist(p, \Delta) + s c_* (v).$$
 which establishes the escape rate inequality (\ref{escape rate})  for the ray  $\gamma(s) = p + s v$ with escape rate 
 $c = c_* (v)$.  Fix any $t >0$.  
 the escape rate  inequality shows that if $p \in N_t (\Delta)$ then the curve $p + s v$  has left $N_t (\Delta)$,
 by the time  $s = t/c_* (v)$,so   that the escape rate from $\Delta$ is positive and at least
 $c_* (v)$.  
 
 To verify the claim regarding the precise global   escape rate observe that $c_* (v) = dist(v, \Delta)$.
 Now take the   maximum of $c_* (v)$ over all  $v \in \S \cap K$.   
Since both $dist(v, \Delta)$ and $\|v \|$ are homogeneous of degree $1$, this  maximum of $c_*(v)$ 
 equals  the maximum of $dist(v, \Delta)/ \|v \|$ over   $v \in K \setminus \{0\}$.
 This last  maximum in turn is  the reciprocal of the {\it  minimum}  of $\|q \|/ dist(q, \Delta)$ over   $q \in K \setminus \{0\}$.
 We can understand this minimum by setting $dist(q, \Delta) = 1$ and then minimizing $\|q \|$.
 But this minimum value is $dist(0, K_1)$
 and is achieved as $\|q_* \|$ where  $q_*$ is the closest point to $0$ on $K_1$. 
 It follows that the    optimal escape rate is equal to  $1/dist(0, K_1)$
 and that the corresponding translational ray escapers $p + s v$
 are obtained by setting  $v = q_*/\|q_* \|$.
 
 To show that this  value just computed  for the  escape rate is best possible 
 consider the problem of escaping from the cone point $p = 0$ into $K$.   The best $t=1$-escaper for $p =0$ will be the shortest 
 path to $\partial K_1$ which is the segment $[0, q_*]$.  Now, by  homogeneity the intersection of the  ray $s  v,  v = q_*/\|q_*\|$
 with $K_t$, $t > 0$,  
 yields the best   t-escaper from $0$,  for all $t$.  
 
 QED

\vskip 2cm

  \section{Appendix: On the distance formula.}
  
  Here we prove the distance formula (\ref{distance fmla}) used in the paper \footnote{
  Doubtless this formula is proved in other papers, probably in some of my own,
  but, not finding a derivation,   I opted to give this one.}.
  
  Fix a point $q = (q_1, q_2, \ldots , q_N)  \in \E$.  The distance between $q$ and any affine subspace $S$ is
  the length of the unique line segment  which hits $S$ orthogonally.
  We apply this observation to $S = \Delta_{ab}$.  A line segment from $q$ to $S$ can
  be written $\ell(t) = (1-t) q + t s$ for some $s \in S$. 
  
  For simplicity of notation, suppose that $a =1, b = 2$. 
  Write 
  $$q_{cm} = \frac{1}{m_1 + m_2}(m_1 q_1 + m_2 q_2) $$
  for the center of mass of 1 and 2.
  The (unique) point of $\Delta_{12}= S$ which we will want with our line segment
  from $q$ turns out to be  
  $$s = (q_{cm}, q_{cm}, q_3, q_4, \ldots, q_N).$$
  The corresponding line  segment is then 
  $\ell (t) = (1-t) q + t s$, or, in terms of components 
  $$\ell_1 (t) = (1-t) q_1 + t q_{cm}$$
  $$\ell_2 (t) = (1-2) q_2 + t q_{cm}$$
  and 
  $$\ell_a (t) = q_a = const. , \text{ for } a > 2.$$
  The corresponding velocity $v = \dot \ell = (v_1, v_2, \ldots, v_N) \in \E$
is given by 
  $$v_1 = -q_1 + q_{cm}$$
  $$v_2 = -q_2 + q_{cm}$$
  and
  $$v_a = 0,  \text{ for } a > 2.$$
  
  The total linear momentum of this trajectory $\ell (t)$ is zero.
Indeed  this total linear momentum, $\Sigma m_a v_a$  is 
   $m_1 v_1 + m_2 v_2 =   -m_1  q_1 - m_2 q_2 + m_{12} q_{cm} = 0$.
   Now let  $h$ be any vector in $\Delta_{12}$.  Then $h$ is of the form 
   $h = (a, a, h_3, h_4, \ldots , h_N)$ where $a, h_i \in \R^d$ are arbitrary.
   We compute that
   $<v, h > = m_1 a v_1 + m_2 a v_2 = a \cdot (m_1 \dot v_1 + m_2 v_2 ) = 0$.
   which shows that the line segment is orthogonal to $\Delta_{12}$. 
   
   We have established that the path $\ell$ is the (unique) line segment joining $q$ to $\Delta_{12}$
 orthogonally.   The length of $\ell$ is thus $dist(q, \Delta_{12})$.
   But this length is  $\|v \|$ where
   $\|v \|^2 =  m_1 |v_1|  + m_2 |v_2|$.   
   The computation is finished with some algebra.
   
   The algebra can be streamlined by introducing the ``probabilities''  $p_1 = \frac{m_1}{m_1 + m_2}$
   and $p_2 = \frac{m_2}{m_1 + m_2}$ which allow us to write
   $$q_{cm} = p_1 q_1 + p_2 q_2,  q_1 = p_1 q_1 + p_2 q_1 ,  q_2 = p_1 q_2 + p_2 q_2.$$
   from which we compute that
   $$v_1 = -p_2 (q_1 - q_2)$$
   $$v_2 = p_1 (q_1 -q_2)$$
   and finally $\|v \|^2 = ( m_1 p_2 ^2 + m_2 p_1 ^2)r_{12} ^2$.
   A wee bit of algebra yields that $( m_1 p_2 ^2 + m_2 p_1 ^2) = \frac{m_1 m_2}{m_1 + m_2}$
   which is the claim.
   
   \section{A cone and its equidistant}

 I was temporarily seduced into the misbelief  that $K_t$, $t > 0$,  must be a cone  since $K_0$ is a cone.
 This is false.  See  figure \ref{fig: coneEg}.   
 
 Correcting my misbelief  corrected my
 intuition and helped me  come up with the  short correct proof given here.   
 For a simple example where $K_1$ is not a cone, suppose that $K_0$ is  the cone over the  
 with sides $2/a$ and $2/b$ where $a, b > 0$ and $a \ne b$.    Place the rectangle on the plane $z = 1$ with
 the cone point at the origin of the xyz plane.  Then we can specify  $K_0$  by the inequalities
 $$z \ge a|x|,  z \ge b |y|.$$
 Its cross sections $z = z_0$  are   rectangles whose sides are in the ratio $1/a : 1/b$
 and which grow linearly with $z_0$.

 To compute $K_t$  find the four normalized  linear functionals which define
 $K_0$.  They are $\ell_{\pm} = \frac{1}{\sqrt{a^2 + 1}} (z \pm a x),  f_{\pm} = \frac{1}{\sqrt{a^2 + 1}}(z \pm b y)$.
 Thus $K_0$ is defined by $\ell_{\pm} \ge 0$ and $f_{\pm} \ge 0$,
 and $K_t$ by $\ell_{\pm} \ge t,  f_{\pm} \ge t$.     Taking $t = 1$
 and doing a bit of algebra yields that $K_1$ is defined by
 $$z - \sqrt{a^2 + 1} \ge a|x|, z  - \sqrt{b^2 + 1} \ge b|y|.$$
 In particular, since the right hand side of these equations is greater than or equal to zero
 we have that $z \ge max \{ \sqrt{a^2 + 1} ,  \sqrt{b^2 + 1}  \}$. 
 
 For concreteness,   set $a = 1$ and suppose   $b < 1$
 so that $z \ge \sqrt{2}$ on $K_1$.  The cross-section $K_1$ with the
 plane $z = \sqrt{2}$ is the bounded  interval $\sqrt{2} - \sqrt{b^2 +1} \ge b |y|$
 or $\frac{1}{b} (\sqrt{2} - \sqrt{b^2 +1} ) \ge   |y|$. 
 Consequently $K_1$ cannot be a cone.  To get a 3-dimensional
 picture of $K_1$ consider the cross-sections $z = z_0$ for  $z_0 > \sqrt{2}$  of $K_1$.
 These are   rectangles $ z_0 - \sqrt{2} \ge |x|$ and 
 $\frac{1}{b} (z_0 - \sqrt{b^2 +1} ) \ge   |y|$.  Working out the side lengths
 we see that the  rectangles have 
 aspect ratio $ (1 - \frac{\sqrt{2}}{z_0}) : ( \frac{1}{b} - \frac{\sqrt{b^2 + 1}}{ z_0 b})$
% which begins at $0:1$ with $z_0 = \sqrt{2}$
 and asymptotes to $1: \frac{1}{b}$, the aspect ratio of the rectangle on which $K_0$ is
 based.   
 
\section{Acknowledgements}   Andrei Agrachev pointed the way towards  reducing the escape game to the  convex case.  Saul Rodriguez Martin, as user Sa\'ul RM  on Math Overflow,    described   the equidistant  $K_1$ of  Appendix B.    I would also like to thank Serge Tabachnikov
 and J.M. Burgos for helpful email conversations.


\begin{thebibliography}{50}
 
     
\bibitem{AbMa} Abraham, R.  and Marsden, J. E.,  (1978)
{\bf Foundations of Mechanics},  2nd ed.   Benjamin/Cummings.


\bibitem{Albouy-Chenciner}  Albouy, A. and Chenciner, A., (1998):  {\em Le probl\'eme des N corps et les distances mutuelles}   Inventiones  {\bf 131} , 151-184
\url{http://link.springer.com/article/10.1007%2Fs002220050200}

\bibitem{Chenciner} Chenciner, A.,  (1997)  {\em \`A l'infini en temps fini}  Seminaire Bourbaki, vol. 1996-7, expos\'es 832, pp 323-353, 
in Ast\'erisque, no. 245, 
\url{http://www.numdam.org/item/?id=SB_1996-1997__39__323_0}

 
 \bibitem{Burago} D. Burago, Y. Burago, S. Ivanov, {\bf A course in metric geometry}  Graduate Studies in Mathematics, {\bf 33} American Mathematical Society, Providence, RI, 2001.
 
 \bibitem{Burgos1}   Burgos, J. M.,: {\em Existence of partially hyperbolic motions in the N-body problem}  Proc. Amer. Math. Soc.  
{\bf 150} , no.  4, April 2022,   1729 --1733
\url{https://doi.org/10.1090/proc/15778}
  
 \bibitem{Burgos2}     Burgos, J. M. and  Maderna, E., {\em Geodesic rays of the  N -body problem}
Arch. Ration. Mech. Anal. {\bf 243}  (2022), no. 2, 807-827.
 
 
\bibitem{Knauf} Knauf, A., (2012)  {\bf Mathematical Physics: Classical Mechanics},
Springer. 


\bibitem{Landau} Landau, I. and Lifshitz, E.,  (1976), {\bf Mechanics}, Pergamon Press.

 \bibitem{Maderna1} Maderna, E., Venturelli, A.: {\em Viscosity solutions and hyperbolic motions: A new PDE method for the N-body problem} Ann. of Math. (2) 192(2), 499--550, 2020 
 
  
 \bibitem{Moeckel1} Moeckel, R.  {\em Minimal geodesics of the isosceles three body problem},
Qual. Theory Dyn. Syst. 19 (2020), no. 1, Paper No. 48, , 29 pp.



 \bibitem{Brake-to-Syz} Moeckel, R. Montgomery, R. and  Venturelli A., {\em From Brake to Syzygy}, Archive for Rational Mechanics and Analysis; {\bf 204}, no. 3, pp. 1009-1060, 2012.
 

 \bibitem{} Montgomery, R., {\em Minimizers for the Kepler problem}, 
Qual. Theory Dyn. Syst. 19 (2020), no. 1, Paper No. 31, 12 pp.

 \bibitem{} Montgomery, R., {\em Who's afraid of the Hill boundary?}, 
SIGMA Symmetry Integrability Geom. Methods Appl. 10 (2014), Paper 101, 11 pp.

 \bibitem{MontRCD}  Montgomery, R., {\em Brake Orbits Fill the N-Body Hill Region}, Regular and Chaotic Dynamics,
{\bf 28} , Nos. 4-5,  374-394, 2023.    Proposition 1. 

\bibitem{Terracini}    Polimeni, C.  and  Terracini, S.,   {\em On the existence of minimal expansive solutions to the N-body problem}
2024 


\bibitem{Rockafellar}  Rockafellar, R. T.,  {\bf Convex Analysis} ,  Princeton U. Press, 1997, 
reprinted from the 1970 edition.
 
 
 \end{thebibliography}
\end{document}